\documentclass[11pt]{amsart}
\usepackage[utf8]{inputenc}
\usepackage{amsmath,amsfonts,amssymb,amsthm}
\usepackage{latexsym,bm,graphicx}
\usepackage{mathrsfs}
\usepackage{color}
\usepackage{hyperref}

\title[Local isoperimetric inequalities]{Local isoperimetric inequalities in metric measure spaces verifying measure contraction property}
\author{Xian-Tao Huang}
\address{School of Mathematics\\  Sun Yat-sen University\\ Guangzhou 510275\\ E-mail address: hxiant@mail2.sysu.edu.cn}

\newtheorem{thm}{Theorem}[section]

\newtheorem{defn}[thm]{Definition}
\newtheorem{rem}[thm]{Remark}

\numberwithin{equation}{section}

\begin{document}

\maketitle
\begin{abstract}
We prove that on an essentially non-branching $\mathrm{MCP}(K,N)$ space, if a geodesic ball has a volume lower bound and satisfies some additional geometric conditions,
then in a smaller geodesic ball (in a quantified sense) we have an estimate on the isoperimetric constants.

\vspace*{5pt}

\vspace*{5pt}

\end{abstract}
\section{Introduction}  

The isoperimetric problem is one of the most classical and beautiful problems in mathematics.
It addresses the following natural problem: given a space $X$, what is the minimal amount of area needed to enclose a fixed volume $v$?

If $X$ is $\mathbb{R}^{N}$, then it is well known that, for every finite perimeter subset $E\subset X$, it holds
$$|\partial E|\geq N\omega_{N}^{\frac{1}{N}}|E|^{\frac{N-1}{N}},$$
(where $|\partial E|$ and $|E|$ denote the $N-1$ and $N$ dimensional volume respectively) and the only optimal shapes
are the round balls.
If $X$ is a manifold with many symmetries such as $S^{N}$ and $\mathbb{H}^{N}$, or is a perturbation of them, there are also plenty of works concerning the isoperimetric problem and describing the `optimal shapes' on it.
The readers can refer to Appendix H in \cite{EM13} for a list of references.

If $X$ is a general manifold, one can only hope some comparison results for the isoperimetric problem (under some curvature assumptions on $X$).
In this direction, the famous L\'{e}vy-Gromov isoperimetric inequality (see Appendix C in \cite{Gro07}) says if $X$ is an $N$-dimensional manifold with Ricci curvature bounded from below by $K>0$, and suppose $E\subset X$ is a finite perimeter subset, then we have
$$\frac{|\partial E|}{|X|}\geq \frac{|\partial B|}{|S|},$$
where $S$ is the $N$-dimensional round sphere with Ricci curvature $K$, and $B\subset S$ is a spherical cap such that $|E|/|X|=|B|/|S|$.
Some extensions of L\'{e}vy-Gromov inequality can be found in \cite{BBG85} \cite{Cr82} etc.
In \cite{Mil15}, E. Milman obtained sharp isoperimetric inequalities which extends the L\'{e}vy-Gromov inequality to smooth manifolds with densities which have generalized Ricci curvature at least $K\in\mathbb{R}$, generalized dimension at most $N\geq 1$ and diameter at most $D<+\infty$.

We note that, given any point on a Riemannian manifold, there is sufficiently small neighborhood around it which is sufficiently close to a ball in Euclidean space, thus in such a neighborhood, the local isoperimetric constant is close to the Euclidean one.
Under some geometric conditions, we may even obtain some quantified estimate in the form that, on a geodesic ball of definite radius, the isoperimetric constant is close to the one in Euclidean ball in a quantified sense.
See Remark 10.5 in \cite{Per02} and Theorem 1.1 in \cite{CaMo17}.
Such an almost-Euclidean isoperimetric inequality is useful in some other problems.
For example, in Perelman's Pseudo-Locality Theorem (see \cite{Per02}), almost-Euclidean isoperimetric inequality is used to obtain curvature estimates for the Ricci flow.

\vspace*{10pt}

One can also consider the isoperimetric problem when $X$ is not a Riemannian manifold.
In fact, some of the above mentioned isoperimetric problem has been considered on non-Riemannian manifolds or even  metric measure spaces.

Recently, people are more and more interested in the study of non-smooth objects, and there are lots of researches on the notion of Ricci curvature lower bounds on metric measure spaces.
Using the theory of optimal transformation, the so-called $\textmd{CD}(K,N)$-condition, which is a notion to describe  `Ricci curvature bounded from below by $K\in \mathbb{R}$ and dimension bounded above by $N\in[1,\infty]$' for general metric measure spaces, was introduced independently by Lott and Villani (\cite{LV09}) and by Sturm (\cite{St06I} \cite{St06II}).
The $\textmd{CD}(K,N)$-condition is compatible with the classical curvature-dimension notions on Riemannian manifolds.
Later on, some variant versions of curvature-dimension condition were introduced by some authors, among them,
the Measure Contraction Property $\mathrm{MCP}(K,N)$ was introduced independently by Ohta in \cite{Oh07} and Sturm in \cite{St06II} as a weaker variant of $\mathrm{CD}(K, N)$.
There are many metric measure spaces verifying $\mathrm{MCP}$ condition but not any $\mathrm{CD}$ condition: e.g. the Heisenberg groups, generalized $H$-type groups, the Grushin plane and Sasakian structures (under some curvature bounds) etc., for more details, see \cite{Jui09} \cite{BaRi18} \cite{BaRi19} \cite{LLZ16} etc.
Thus, researches on general $\mathrm{MCP}(K,N)$ spaces may give information which are new even on the above mentioned examples.

In \cite{CaMo17I}, Cavalletti and Mondino extended the L\'{e}vy-Gromov-Milman isoperimetric inequality to the class of essentially non-branching (see Section \ref{sec2} for the definition) metric measure spaces verifying $\mathrm{CD}(K, N)$ with $m(X)=1$.
The key tool in \cite{CaMo17I} is the localization technique, which is mainly based on the work developed by Payne-Weinberger \cite{PaWe60}, Gromov-Milman \cite{GrMi87}, Kannan-Lov\'{a}sz-Simonovits \cite{KLS87} and Klartag \cite{Kla17}.
In a word, in \cite{CaMo17I}, using the theory of $L^{1}$-Optimal Transport, the authors transform the isoperimetric problem on a $\mathrm{CD}(K,N)$ space to the isoperimetric problem on one-dimensional $\mathrm{CD}(K,N)$ spaces.

Using the localization technique again, Cavalletti and Mondino studied the local isoperimetric inequality in essentially non-branching $\mathrm{CD}(K,N)$ spaces in \cite{CaMo17} and obtain the following theorem:

\begin{thm}[Theorem 1.4 in \cite{CaMo17}]\label{thmCM1}
For every $K\in\mathbb{R}, N\in [2,\infty)$ there exist $\bar{\epsilon}_{K,N}$, $\bar{\eta}_{K,N}$,  $\bar{\delta}_{K,N}$, $C_{K,N} > 0$ such that the next statement is satisfied.
Let $(X, d, m)$ be a geodesic metric space endowed with a non-negative Borel measure.
For a fixed $\bar{x}\in X$, assume that $B_{4\bar{r}_{K,N}}(\bar{x})$ is relatively compact and that $B_{4\bar{r}_{K,N}}(\bar{x})\subset \mathrm{supp}(m)$, $m(B_{4\bar{r}_{K,N}}(\bar{x}))<\infty$.
Assume moreover that for some $\epsilon\in[0, \bar{\epsilon}_{K,N}],\eta\in[0, \bar{\eta}_{K,N}]$ it holds:
\begin{description}
  \item[(1)] $m(B_{\bar{r}_{K,N}}(\bar{x}))\geq 1-\eta$;
  \item[(2)] $\limsup_{r\downarrow0}\frac{m(B_{r}(\bar{x}))}{\omega_{N}r^{N}}\leq1+\eta$;
  \item[(3)] $(X,d,m)$ is essentially non-branching and verifies $\mathrm{CD}_{\mathrm{loc}}(K-\epsilon,N)$ inside $B_{4\bar{r}_{K,N}}(\bar{x})$.
\end{description}
Then for every $\delta\in(0, \bar{\delta}_{K,N}]$ and every finite perimeter subset $E\subset B_{\delta}(\bar{x})$ the following almost Euclidean isoperimetric inequality holds:
\begin{align}\label{1.1}
P(E)\geq N\omega_{N}^{\frac{1}{N}}(1-C_{K,N}(\delta+\epsilon+\eta))m(E)^{\frac{N-1}{N}}.
\end{align}
\end{thm}

For $N\in(1, 2)$, as is pointed out in Remark 1.5 of \cite{CaMo17}, a conclusion similar to (\ref{1.1}) holds with a bit difference: the power on $\delta$ in the error term is $\frac{2(N-1)}{N}$ in this case.

Note that Theorem \ref{thmCM1} recovers a theorem claimed by Perelman in \cite{Per02}, see Theorem 1.1 in \cite{CaMo17}.

In the following, we explain some notation appeared in the statement of Theorem \ref{thmCM1}, some similar notations also appear in the rest part of the paper.

We say $(X,d,m)$ verifies $\mathrm{CD}_{\mathrm{loc}}(K-\epsilon,N)$ inside $B_{4\bar{r}_{K,N}}(\bar{x})$ if for every $x \in B_{4\bar{r}_{K,N}}(\bar{x})$, there exists a neighbourhood $U$ such that $\mathrm{CD}(K-\epsilon,N)$ is verified inside $U$ (see \cite{BS10}).

For any $N\in(1,\infty)$, we define the function $r\mapsto \mathrm{Vol}_{K,N}(r)$ to be:

\begin{align}\label{1.2}
\mathrm{Vol}_{K,N}(r):=
   \left\{
     \begin{array}{ll}
       N\omega_{N}\int_{0}^{r}\sin(t\sqrt{\frac{K}{N-1}})^{N-1}dt, & \hbox{if $K>0$;} \\
       \omega_{N}r^{N}, & \hbox{if $K=0$;} \\
       N\omega_{N}\int_{0}^{r}\sinh(t\sqrt{\frac{K}{N-1}})^{N-1}dt, & \hbox{if $K<0$.}
     \end{array}
   \right.
\end{align}
where
$$\omega_{N}:=\frac{\pi^{\frac{N}{2}}}{\Gamma(\frac{N}{2}+1)},$$
with $\Gamma$ denoting the Euler's Gamma function.
If $N\in \mathbb{N}$ then $\mathrm{Vol}_{K,N}(r)$ is nothing but the volume of the metric ball of radius $r$ in $\mathbb{M}^{N}_{K/(N-1)}$, the simply connected manifold of constant sectional curvature equal to $\frac{K}{N-1}$.

Then the positive constant $\bar{r}_{K,N}$ is defined so that $\mathrm{Vol}_{K,N}(\bar{r}_{K,N})=1$.

For a subset $E\subset X$, the perimeter of $E$ is defined to be
$$P(E)=\inf\bigl\{\liminf_{h\rightarrow\infty}\int_{X}\mathrm{lip}(f_{h})dm \mid f_{h}\in\mathrm{Lip}(X), \lim_{h\rightarrow0}\int_{X}|f_{h}-\chi_{A}|dm=0\bigr\},$$
where $\mathrm{lip}(f_{h})$ is the local Lipschitz constant for a Lipschitz function $f_{h}$.
\vspace*{10pt}

In \cite{CaSa18}, Cavalletti and Santarcangelo considered isoperimetric inequalities on $\mathrm{MCP}(K,N)$ spaces, and they have obtained sharp L\'{e}vy-Gromov type isoperimetric inequalities on essentially non-branching $\mathrm{MCP}(K,N)$ spaces with diameter upper bound $D$.
In fact, the localization technique still applies to essentially non-branching $\mathrm{MCP}(K,N)$ spaces, see e.g. \cite{BC09}, \cite{CaMi16}, \cite{CaMo18} etc.
Using localization, the isoperimetric problem on a $\mathrm{MCP}(K,N)$ space is transformed to the corresponding statement on one-dimensional $\mathrm{MCP}(K,N)$ spaces.
In \cite{CaSa18}, the authors obtain the explicit description of the optimal one-dimensional $\mathrm{MCP}(K, N)$-density and study some fundamental properties of these densities.

Recently, there are many other researches on $\mathrm{MCP}(K,N)$ spaces basing on localization technique, see \cite{H19} \cite{HaMi18} etc.

\vspace*{10pt}

Motivated by \cite{CaMo17}, in this paper, we consider the local isoperimetric constant for essentially non-branching $\mathrm{MCP}(K,N)$ spaces.

From now on, $(X, d, m)$ will be an essentially non-branching $\mathrm{MCP}(K,N)$ space with $\mathrm{supp}(m)=X$, where $N>1$.
Given any fixed $D>0$ and $x\in X$, we define a function
$$f^{X}_{x,D}(r)=\frac{1}{m(B_{D}(x))}m(B_{r}(x)).$$

We use the notation $\Psi(u_{1},\ldots,u_{k}\mid \ldots)$ to denote a nonnegative function depending on the
numbers $u_{1},\ldots,u_{k}$ and some additional parameters, such that when these additional parameters are fixed, we have
$$\lim_{u_{1},\ldots,u_{k}\rightarrow0}\Psi(u_{1},\ldots,u_{k}\mid \ldots)=0.$$

The following two theorems are main results of this paper:

\begin{thm}\label{thm-local-isoperi1}
Given $N>1$, and let $K=0$ or $-(N-1)$.
Fix $D>0$ and a function $\bar{f}:(0,D)\rightarrow\mathbb{R}^{+}$ with $\lim_{r\downarrow0}\bar{f}(r)=0$.
There exists $\bar{\delta}> 0$ depending on $N,D,\bar{f}$ such that the next statement is satisfied.
Suppose $(X, d, m)$ is an essentially non-branching $\mathrm{MCP}(K,N)$ space,
$\bar{x}\in X$.
Assume in addition:
\begin{description}
  \item[(1)] $m(B_{D}(\bar{x}))\geq \mathrm{Vol}_{K,N}(D)$;
  \item[(2)] there exists $r_{0}>0$ such that $f^{X}_{\bar{x},D}(r)\leq \bar{f}(r)$ for every $r\in(0, r_{0})$.
\end{description}
Then for every $\delta\in(0, \bar{\delta}]\cap (0, r_{0}]$ and every finite perimeter subset $E \subset B_{\delta}(\bar{x})$, the following isoperimetric inequality holds:
\begin{align}\label{1.3}
P(E)\geq (1-\Psi(\delta\mid N, D, \bar{f}))N^{\frac{1}{N}}\omega_{N}^{\frac{1}{N}}m(E)^{\frac{N-1}{N}}.
\end{align}
\end{thm}

\begin{thm}\label{thm-local-isoperi2}
Given $N>1$, and let $K=N-1$.
Fix $D\in(0,\pi)$.
There exist $\bar{\eta}, \bar{\delta}> 0$ depending on $N,D$ such that the next statement is satisfied.
Suppose $(X, d, m)$ is an essentially non-branching $\mathrm{MCP}(K,N)$ space,
$\bar{x}\in X$.
Assume moreover that for some $\eta\in[0, \bar{\eta}]$ it holds:
\begin{description}
  \item[(1)] $m(B_{D}(\bar{x}))\geq \mathrm{Vol}_{K,N}(D)$;
  \item[(2)] $\limsup_{r\downarrow0}\frac{m(B_{r}(\bar{x}))}{\omega_{N}r^{N}}\leq1+\eta$.
\end{description}
Then for every $\delta\in(0, \bar{\delta}]$ and every finite perimeter subset $E\subset B_{\delta}(\bar{x})$, the following isoperimetric inequality holds:
\begin{align}\label{1.4}
P(E)\geq (1-\Psi(\delta,\eta\mid N, D))N^{\frac{1}{N}}\omega_{N}^{\frac{1}{N}}m(E)^{\frac{N-1}{N}}.
\end{align}

\end{thm}

\begin{rem}
\begin{description}
  \item[(1)] In the statement of Theorem \ref{thmCM1}, $K$ is a variable quantity and $\bar{r}_{K,N}$ is chosen so that $\mathrm{Vol}_{K,N}(\bar{r}_{K,N})=1$; while in the statement of Theorem \ref{thm-local-isoperi1} and Theorem \ref{thm-local-isoperi2}, $K$ is normalized but $D$ is variable.
      The two statements are equivalent as we can consider the rescaling $(X, \lambda_{1} d, \lambda_{2}m)$ for suitable $\lambda_{1}$ and $\lambda_{2}$.
  \item[(2)] There is a bit difference in assumption (1) of Theorem \ref{thm-local-isoperi1} and Theorem \ref{thm-local-isoperi2} and assumption (1) of Theorem \ref{thmCM1}, but the two statements are equivalent.
      In fact, if assumption (1) of Theorem \ref{thm-local-isoperi1} or Theorem \ref{thm-local-isoperi2} is replaced by $m(B_{D}(\bar{x}))\geq (1-\eta)\mathrm{Vol}_{K,N}(D)$, then we can consider $(X,d,\tilde{m})=(X,d,\frac{1}{1-\eta}m)$.
      It is easy to see, in the case $K=-(N-1)$ or $0$, $(X,d,\tilde{m})$ satisfies all the assumptions in Theorem \ref{thm-local-isoperi1}, hence
      \begin{align}\label{1.5}
      P(E)&=(1-\eta)\tilde{P}(E)\geq (1-\eta)(1-\Psi(\delta\mid N, D,\bar{f}))N^{\frac{1}{N}}\omega_{N}^{\frac{1}{N}}\tilde{m}(E)^{\frac{N-1}{N}}\\
      &= (1-\eta)^{\frac{1}{N}}(1-\Psi(\delta\mid N, D,\bar{f}))N^{\frac{1}{N}}\omega_{N}^{\frac{1}{N}}m(E)^{\frac{N-1}{N}},\nonumber
      \end{align}
      where $\tilde{P}(E)$ is the perimeter of $E$ in $(X,d,\tilde{m})$.
      In the case $K=N-1$, $(X,d,\tilde{m})$ satisfies $\limsup_{r\downarrow0}\frac{\tilde{m}(B_{r}(\bar{x}))}{\omega_{N}r^{N}}\leq\frac{1+\eta}{1-\eta}\leq 1+3\eta$ provided $\eta$ is sufficiently small, then we apply Theorem \ref{thm-local-isoperi2} to $(X,d,\tilde{m})$ and obtain
      \begin{align}
      P(E)\geq (1-\Psi(\delta,\eta\mid N, D))N^{\frac{1}{N}}\omega_{N}^{\frac{1}{N}}m(E)^{\frac{N-1}{N}} \nonumber
      \end{align}
      similar to (\ref{1.5}).
      On the other hand, if we have proved theorems with assumption (1) replaced by $m(B_{D}(\bar{x}))\geq (1-\eta)\mathrm{Vol}_{K,N}(D)$, then it is easy to see Theorem \ref{thm-local-isoperi1} and Theorem \ref{thm-local-isoperi2} also hold.
  \item[(3)] In (\ref{1.3}), the principal coefficient $N^{\frac{1}{N}}\omega_{N}^{\frac{1}{N}}$ is smaller than the one in manifolds (or in $\mathrm{CD}(K,N)$ spaces as in Theorem \ref{thmCM1}).
      But this constant is sharp in the class of $\mathrm{MCP}(K,N)$ spaces (for the $K\leq0$ cases), as it is almost attained by a class of $1$-dimensional $\mathrm{MCP}(K,N)$ spaces,
      see Remark \ref{rem5.2}.
      In Theorem \ref{thm-local-isoperi2}, we assume (2) because of technical reasons, but (2) is not satisfied by the $1$-dimensional $\mathrm{MCP}(K,N)$ spaces in Remark \ref{rem5.2}.
      It may be interesting to drop assumption (2) in Theorem \ref{thm-local-isoperi2} or to improve the principal coefficient in (\ref{1.4}) under assumption (2).
\end{description}
\end{rem}

\vspace*{10pt}

In the proofs of Theorem \ref{thm-local-isoperi1} and Theorem \ref{thm-local-isoperi2}, we will only handle the isoperimetric inequalities for the outer Minkowski content.
Here we recall that, for a subset $E\subset X$, its outer Minkowski content is defined to be
$$m^{+}(E)=\liminf_{\epsilon\rightarrow0}\frac{m(E^{\epsilon})-m(E)}{\epsilon},$$
where $E^{\epsilon}:=\{x\in X|d(x, E) <\epsilon\}$.
In fact, it is proved in \cite{AGM16} that, on general metric measure spaces, the perimeter is the relaxation of the outer Minkowski content with respect to convergence in measure, hence isoperimetric inequalities written in terms of the outer Minkowski content are equivalent to the corresponding statements written in terms of the perimeter.
See also \cite{CaMo18I} for related results under curvature assumptions.

Our proofs of Theorem \ref{thm-local-isoperi1} and Theorem \ref{thm-local-isoperi2} are based on the localization technique on essentially non-branching $\mathrm{MCP}(K,N)$ spaces.
Similar to Cavalletti and Mondino's paper \cite{CaMo17}, there are some necessary modifications when applying the localization technique.
By the localization technique, the local isoperimetric problem is reduced to some analytic problems on $1$-dimensional $\mathrm{MCP}(K,N)$ spaces.
Note that in \cite{CaSa18}, the optimal shapes for the isoperimetric problem for $1$-dimensional $\mathrm{MCP}(K,N)$ spaces are studied.
In order to handle the local isoperimetric problem, we need to obtain some new properties on the $1$-dimensional spaces, see Section \ref{sec3}.

\vspace*{10pt}

\noindent\textbf{Acknowledgments.} The author would like to thank Prof. X.-P. Zhu, H.-C. Zhang and B.-X. Han for discussions.
The author is partially supported by NSFC (Nos. 12025109 and 11521101) and Guang-dong Natural Science Foundation 2019A1515011804.

\section{Preliminaries}\label{sec2}

Throughout this paper, we will always assume the metric measure space $(X,d,m)$ we consider satisfies the following: $(X, d)$ is a complete separable locally compact geodesic space, and $m$ is a nonnegative Radon measure with respect to $d$ and finite on bounded sets, $\textmd{supp}(m)=X$.

A curve $\gamma:[0,T]\rightarrow X$ is called a geodesic provided $d(\gamma_{s},\gamma_{t})=L(\gamma|_{[s,t]})$ for every $[s,t]\subset[0,T]$, where $L(\gamma)$ means the length of the curve $\gamma$.
$(X,d)$ is called a geodesic space if every two points $x,y\in X$ are connected by a geodesic $\gamma$.
$\textmd{Geo}(X)$ denotes the set of all geodesics with domain $[0,1]$.
For $t\in[0,1]$, define the evaluation map $e_{t}:\textmd{Geo}(X)\rightarrow X$ by $e_{t}(\gamma)=\gamma_{t}$.

Denote by $\mathcal{P}(X)$ the space of Borel probability measures on $X$, and $\mathcal{P}_{2}(X)\subset \mathcal{P}(X)$ the space of Borel probability measures $\xi$ satisfying $\int_{X} d^{2}(x,y)\xi(dy)<\infty$ for some (and hence all) $x\in X$.

For $\mu,\nu\in{\mathcal{P}_{2}(X)}$, consider their Wasserstein distance $W_{2}(\mu,\nu)$ defined by
\begin{align}\label{2.4}
W_{2}^{2}(\mu,\nu)=\underset{\eta\in\Pi(\mu,\nu)}{\inf}\int_{X\times{X}}d^{2}(x,y)d\eta(x,y),
\end{align}
where $\Pi(\mu,\nu)$ is the set of Borel probability measures $\eta$ on $X\times{X}$ satisfying $\eta(A\times{X})=\mu(A)$, $\eta(X\times{A})=\nu(A)$ for every Borel set $A\subset{X}$.
It is known that the infimum in (\ref{2.4}) is always attained for any $\mu,\nu\in{\mathcal{P}_{2}(X)}$.
See \cite{AG11} \cite{Vi09} for the theory of optimal transport.

Given $\mu_{0}, \mu_{1}\in\mathcal{P}_{2}(X)$, we denote by $\mathrm{OptGeo}(\mu_{0}, \mu_{1})$ the space of all $\Pi\in\mathcal{P}(\mathrm{Geo}(X))$ for which $(e_{0}, e_{1})_{\#}\Pi$ realizes the minimum in (\ref{2.4}).
If $(X, d)$ is geodesic, then the set $\mathrm{OptGeo}(\mu_{0}, \mu_{1})$ is non-empty for any $\mu_{0}, \mu_{1}\in\mathcal{P}_{2}(X)$.

We say a subset $G\subset\mathrm{Geo}(X)$ is non-branching if any $\gamma^{1},\gamma^{2}\in G$ with  $\gamma^{1}|_{I}=\gamma^{2}|_{I}$ for some interval $I\subset[0,1]$ must satisfy $\gamma^{1}\equiv\gamma^{2}$ on $[0,1]$.

\begin{defn}
$(X, d, m)$ is called essentially non-branching if and only if for any $\mu_{0}, \mu_{1}\in\mathcal{P}_{2}(X)$, with $\mu_{0}, \mu_{1}\ll m$, any element of $\mathrm{OptGeo}(\mu_{0}, \mu_{1})$ is concentrated on a non-branching subset of geodesics.
\end{defn}

For $\kappa\in \mathbb{R}$, we define the function $s_{\kappa}:[0,+\infty)\rightarrow\mathbb{R}$ (on $[0,\pi/\sqrt{\kappa})$ if $\kappa> 0$) to be
\begin{align}
s_{\kappa}(\theta):=
   \left\{
     \begin{array}{ll}
       (1/\sqrt{\kappa})\sin(\sqrt{\kappa}\theta), & \hbox{if $\kappa>0$;} \\
       \theta, & \hbox{if $\kappa=0$;} \\
       (1/\sqrt{-\kappa})\sinh(\sqrt{-\kappa}\theta), & \hbox{if $\kappa<0$.} \\
     \end{array}
   \right.
\end{align}

Given two numbers $K,N\in\mathbb{R}$ with $N> 1$, for $(t,\theta)\in[0, 1]\times\mathbb{R}^{+}$, we set
\begin{align}
\sigma^{(t)}_{K,N-1}(\theta)
=\left\{
   \begin{array}{ll}
     +\infty, & \hbox{if $K\theta^{2}\geq (N-1)\pi^{2}$;} \\
     \frac{s_{K/(N-1)}(t\theta)}{s_{K/(N-1)}(\theta)}, & \hbox{otherwise.}
   \end{array}
 \right.
\end{align}
and
\begin{align}
\tau^{(t)}_{K,N}(\theta) =t^{\frac{1}{N}}(\sigma^{(t)}_{K,N-1}(\theta))^{\frac{N-1}{N}}.
\end{align}

\begin{defn}\label{def2.8}
We say $(X, d, m)$ satisfies the $(K, N)$-measure contraction property ($\mathrm{MCP}(K, N)$) if for any $x\in X$ and $m$-measurable set $A\subset X$ with $m(A)\in(0,\infty)$,
there exists $\Pi\in\mathrm{OptGeo}(\frac{1}{m(A)}m|_{A},\delta_{x})$, such that
for every $t\in[0,1]$,
  \begin{align}\label{2.10}
  \frac{1}{m(A)}m\geq(e_{t})_{\#}\biggl( \bigl(\tau_{K,N}^{(1-t)}(d(\gamma_{0},x))\bigr)^{N}\Pi(d\gamma)\biggr).
  \end{align}
\end{defn}

\section{Analysis on 1-dimensional model}\label{sec3}

In this section, we consider the isoperimetric problem on $1$-dimensional spaces $(X, d, m) =(I, |\cdot|, h\mathcal{L}^{1})$.

It is well known that $(I, |\cdot|, h\mathcal{L}^{1})$ verifies $\textmd{MCP}(K, N)$ if and only if up to modification on a null-set, the non-negative Borel function $h$ satisfies
$$h(tx_{1} + (1-t)x_{0})\geq\sigma^{(1-t)}_{K,N-1}(|x_{1}-x_{0}|)^{N-1}h(x_{0})$$
for all $x_{0}, x_{1}\in I$ and $t\in[0, 1]$.
We will call $h$ an $\textmd{MCP}(K,N)$ density.
Without loss of generality, we can assume $h$ to be defined over $[0,D]$ for $D\in(0,+\infty]$, and we always assume an $\textmd{MCP}(K,N)$ density $h$ is the continuous representative in its a.e. class (in fact, as a consequence of (2.5) in \cite{CaSa18}, $h$ is locally Lipschitz in the interior of $I$).

Denote by
$$\tilde{\mathcal{F}}_{K,N,D}:=\{\mu\in\mathcal{P}(\mathbb{R})\mid \mathrm{supp}(\mu)=[0,D], \mu= h_{\mu}\mathcal{L}^{1},
h_{\mu} \text{ is an $\mathrm{MCP}(K,N)$ density}\}.$$

For each $v\in(0, 1)$, denote by
$$\tilde{\mathcal{I}}_{K,N,D}(v) :=\inf\{\mu^{+}(A)\mid A\subset[0,D], \mu(A) =v, \mu\in\tilde{\mathcal{F}}_{K,N,D} \}.$$

In \cite{CaSa18}, in order to characterize the optimal shapes for the isoperimetric problem for $1$-dimensional $\mathrm{MCP}(K,N)$ spaces with diameter $\leq D$, the authors define a positive function as follows:
\begin{align}
f_{K,N,D}(x):=\biggl( \int_{0}^{x} \biggl(\frac{s_{K/(N-1)}(D-y)}{s_{K/(N-1)}(D-x)} \biggr)^{N-1}dy +\int_{x}^{D} \biggl(\frac{s_{K/(N-1)}(y)}{s_{K/(N-1)}(x)} \biggr)^{N-1}dy \biggr)^{-1}
\end{align}
if $x\in(0, D)$ and equal $0$ if $x =0, D$.

Then, for each $a\in(0,D)$, let
\begin{align}\label{4.3}
h^{a}_{K,N,D}(x):=f_{K,N,D}(a)
   \left\{
     \begin{array}{ll}
       \bigl(\frac{s_{K/(N-1)}(D-x)}{s_{K/(N-1)}(D-a)} \bigr)^{N-1}, & \hbox{if $x\in[0, a]$;} \\
       \bigl(\frac{s_{K/(N-1)}(x)}{s_{K/(N-1)}(a)} \bigr)^{N-1}, & \hbox{if $x\in [a,D]$.}
     \end{array}
   \right.
\end{align}
In the rest of this paper, the dependence of $h^{a}_{K,N,D}$ on $K, N, D$ will be omitted and we will use $h_{a}$ for simplicity.

One can check that, for each $a\in(0, D)$, $h_{a}$ integrates to $1$ and it is an $\mathrm{MCP}(K, N)$ density, but it does not verify $\mathrm{CD}(K, N)$ condition except the case in which $K>0$ and $D=\pi\sqrt{(N-1)/K}$.
See Lemmas 3.3 and 3.4 in \cite{CaSa18}.

Following \cite{CaSa18}, consider the map
$$(0,D)\ni a\mapsto v_{K,N,D}(a):=\int_{0}^{a}h_{a}(x)dx
\in(0, 1).$$
By Lemma 3.5 in \cite{CaSa18}, $v_{K,N,D}(a)$ is invertible, hence for each $K, N, D$ it is possible to define the inverse map of $v_{K,N,D}(a)$:
$$(0, 1)\ni v\mapsto a_{K,N,D}(v)\in(0,D),$$
with $a_{K,N,D}(v)$ the unique element such that
\begin{align}
\int_{0}^{a_{K,N,D}(v)}h_{a_{K,N,D}(v)}(x)dx=v.
\end{align}

In \cite{CaSa18}, the following theorem is proved:

\begin{thm}[Theorem 3.7 in \cite{CaSa18}]\label{thm4.1}
For each volume $v\in(0,1)$, it holds
\begin{align}\label{4.1}
\tilde{\mathcal{I}}_{K,N,D}(v) = f_{K,N,D}(a_{K,N,D}(v))=h_{a_{K,N,D}(v)}(a_{K,N,D}(v)).
\end{align}
In particular, the lower bound in the definition of $\tilde{\mathcal{I}}_{K,N,D}(v)$ is attained.
\end{thm}

In the following, we will compute $\tilde{\mathcal{I}}_{K,N,D}(v)$ for $v$ sufficiently small.

We will only consider the cases $K=N-1, 0$ and $-(N-1)$, and correspondingly, $\kappa=\frac{K}{N-1}$ take values $1, 0, -1$.
The conclusions for general $K$ can be obtained by rescaling.
We will fix $L\in(0,+\infty)$ and $\lambda\in(0,1]$.
In the $K=N-1$ case, we assume $L<\pi$.
Suppose $D\in[\lambda L, L]$.
Denote by
\begin{align}
k_{D}=\int_{0}^{D}s_{\kappa}(t)dt.
\end{align}

We first fix some notations.
Suppose $g_{1}:(0,\theta_{0})\mapsto \mathbb{R}\setminus\{0\}$ is a function (where $\theta_{0}$ is some positive number), we use $o(g_{1})$ to denote some function $g_{2}:(0,\theta_{0})\mapsto \mathbb{R}$ such that, for any $\epsilon>0$, there exists $\delta_{0}$ depending on $\epsilon$, $N$, the lower and upper bound of $s_{\kappa}$ on $[\lambda L,L]$ and the upper bound of higher order derivatives of $s_{\kappa}$ on $[\lambda L,L]$ (hence depending on $\epsilon$, $N$, $\lambda$ and $L$) such that $|\frac{g_{2}(t)}{g_{1}(t)}|<\epsilon$ holds for every $t\in(0,\delta_{0})$.
In the following, the function $o(g_{1})$ may vary in different lines, but it always satisfies the above mentioned property.

For simplicity, we denote $f_{K,N,D}$, $v_{K,N,D}$ and $a_{K,N,D}$ by $f_{D}$, $v_{D}$ and $a_{D}$ respectively.
By the definition of $f_{D}$, for $x\ll 1$, we have
\begin{align}\label{4.5}
(f_{D}(x))^{-1}&=\frac{\int_{0}^{x} s_{\kappa}(D-y)^{N-1}dy}{s_{\kappa}(D-x)^{N-1}}+\frac{\int_{x}^{D} s_{\kappa}(y) ^{N-1}dy}{s_{\kappa}(x)^{N-1}}  \\
&=\frac{x(s_{\kappa}(D) ^{N-1}+o(1))}{s_{\kappa}(D)^{N-1}+o(1)} + \frac{k_{D}-\int_{0}^{x} s_{\kappa}(y) ^{N-1}dy}{x^{N-1}+o(x^{N-1})} \nonumber \\
&=\frac{x(1+\frac{1}{s_{\kappa}(D) ^{N-1}}o(x))}{1+\frac{1}{s_{\kappa}(D) ^{N-1}}o(x)} + \frac{k_{D}-\int_{0}^{x} (y^{N-1}+o(y^{N-1}))dy}{x^{N-1}+o(x^{N-1})} \nonumber \\
&=x(1+o(x))+ \frac{k_{D}-\frac{1}{N}x^{N}+o(x^{N})}{x^{N-1}+o(x^{N-1})} \nonumber \\
&=x+ \frac{k_{D}}{x^{N-1}+o(x^{N-1})} -\frac{1}{N}x+o(x) \nonumber \\
&=\frac{k_{D}}{x^{N-1}}(1+o(1)), \nonumber
\end{align}
hence
\begin{align}\label{4.10}
f_{D}(x)=\frac{x^{N-1}}{k_{D}}(1+o(1)).
\end{align}

By (\ref{4.3}) and the definition of $v_{D}(a)$, we have
$$v_{D}(a)= \frac{f_{D}(a)}{s_{\kappa}(D-a)^{N-1}}\int_{0}^{a} s_{\kappa}(D-x)^{N-1}dx.$$
If $a\ll1$, then by (\ref{4.10}), we have
\begin{align}\label{4.14}
v_{D}(a)&= \frac{\frac{a^{N-1}}{k_{D}}(1+o(1))}{s_{\kappa}(D)^{N-1}+o(1)}\int_{0}^{a} (s_{\kappa}(D)^{N-1}+o(1))dx\\
&= \frac{a^{N-1}(1+o(1))}{k_{D}s_{\kappa}(D)^{N-1}}a (s_{\kappa}(D)^{N-1}+o(1))\nonumber\\
&= \frac{a^{N}(1+o(1))}{k_{D}}.\nonumber
\end{align}

Note that $v_{D}(a)\rightarrow0$ if $a\rightarrow0$.
Recall that $v_{D}(a)$ is an increasing function (the proof of this fact can be found in the proof of Lemma 3.5 in \cite{CaSa18}).
Hence together with (\ref{4.14}), we can check that $a_{D}(v)\rightarrow0$ if $v\rightarrow0$, and if $v\ll1$,
\begin{align}\label{4.11}
a_{D}(v)=k_{D}^{\frac{1}{N}}v^{\frac{1}{N}}(1+o(1)).
\end{align}

Combining (\ref{4.1}) (\ref{4.10}) and (\ref{4.11}),
\begin{align}\label{4.12}
\tilde{\mathcal{I}}_{K,N,D}(v) = f_{D}(a_{D}(v)) =\frac{(k_{D}^{\frac{1}{N}}v^{\frac{1}{N}})^{N-1}}{k_{D}}(1+o(1)) =k_{D}^{-\frac{1}{N}}v^{\frac{N-1}{N}}(1+o(1))
\end{align}
holds for $v\in (0,\bar{v})$, where $\bar{v}$ is a sufficiently small positive number depending on $N$, $L$ and $\lambda$.

\begin{rem}
We emphasize that in the above argument, we assume $L<\pi$ in the case $K=N-1$.
In \cite{CaSa18}, the definition of $f_{K,N,D}$, $h^{a}_{K,N,D}$, $v_{K,N,D}$ and $a_{K,N,D}$ still make sense in the case $D=\pi$ and $K=N-1$.
But in this case,
\begin{align}
(f_{\pi}(x))^{-1}=\frac{\int_{0}^{x} (\sin(\pi-y))^{N-1}dy}{(\sin(\pi-x))^{N-1}}+\frac{\int_{x}^{\pi} (\sin y) ^{N-1}dy}{(\sin x)^{N-1}} =\frac{\int_{0}^{\pi} (\sin y) ^{N-1}dy}{(\sin x)^{N-1}},
\end{align}
and similar to the above argument, we can prove that,
\begin{align}
\tilde{\mathcal{I}}_{N-1,N,\pi}(v) = f_{\pi}(a_{\pi}(v)) =k_{\pi}^{-\frac{1}{N}}N^{\frac{N-1}{N}}v^{\frac{N-1}{N}}(1+o(1))
\end{align}
holds for $v\in (0,\bar{v})$, where $\bar{v}$ is a sufficiently small constant depending on $N$.
\end{rem}

\section{The localization technique on MCP spaces}\label{sec7}

The proofs of Theorem \ref{thm-local-isoperi1} and Theorem \ref{thm-local-isoperi2} are mainly based on the localization technique on essentially non-branching $\mathrm{MCP}(K,N)$ spaces.
The readers can refer to Section 3 in \cite{CaMo18} for details, and consult \cite{CaMo17I} \cite{Ca18} \cite{Ca14-2} etc. for some related details on $\mathrm{CD}(K,N)$ spaces.
As we are considering the local isoperimetric problem in this paper, there are some necessary modifications when the localization technique are applied to, similar to what Cavalletti and Mondino have done in $\mathrm{CD}(K,N)$ spaces (see \cite{CaMo17}).
Our proof follows the ideas in \cite{CaMo17} closely.
For completeness of exposition, in this section we describe some notations in the construction briefly.
We report the main conclusions of localization technique when modified in our setting, while most of their proofs are omitted except necessary; the readers can refer to \cite{CaMo17} \cite{CaMo18} for the missing details.

In this section, $(X, d, m)$ is an essentially non-branching $\mathrm{MCP}(K,N)$ space ($K\in\mathbb{R}$, $N>1$) with $\mathrm{supp}(m)=X$, $\bar{x}\in X$.
$D>0$ is fixed, and we assume $\delta>0$ is sufficiently small (depending on $K$, $N$ and $D$).

Denote by $\bar{m}=\frac{1}{m(B_{D+2\delta}(\bar{x}))}m\mid_{B_{D+2\delta}(\bar{x})}$.

Given any Borel subset $E\subset B_{\delta}(\bar{x})$ with $\bar{m}(E)>0$, considered the function $f_{E}:X\rightarrow\mathbb{R}$ defined by
\begin{align}\label{5.1}
f_{E}(x):=\chi_{E}(x)-\frac{\bar{m}(E)}{\bar{m}(B_{D}(\bar{x}))}\chi_{B_{D}(\bar{x})}(x).
\end{align}
Obviously $\int f_{E}\bar{m} = 0$.
Denote by $f^{+}_{E}:=\max\{f_{E},0\}$, $f^{-}_{E}:=\max\{-f_{E},0\}$, and
$$c_{E}:=\int f^{+}_{E}\bar{m} =\int f^{-}_{E}\bar{m}>0.$$
Set $\mu_{0}:=\frac{1}{c_{E}}f_{E}^{+}\bar{m}\in\mathcal{P}(X)$, $\mu_{1}:=\frac{1}{c_{E}}f_{E}^{-}\bar{m}\in\mathcal{P}(X)$.
Obviously, $\mu_{0}(E) = \mu_{1}(B_{D}(\bar{x})\setminus E)= 1$.

Consider the $L^{1}$-optimal transportation problem from $\mu_{0}$ to $\mu_{1}$. 
By Kantorovich duality (see Theorem 5.10 in \cite{Vi09}), there exists a $1$-Lipschitz function $\varphi: X\rightarrow \mathbb{R}$, called a Kantorovich potential, such that for any optimal plan $\bar{\pi}\in\Pi(\mu_{0},\mu_{1})$, we have $\bar{\pi}(\Gamma_{0})=1$, where
\begin{align}
\Gamma_{0}:=\{(x, y)\in X\times X\mid \varphi(x)-\varphi(y)=d(x,y)\}.
\end{align}
Define $\Gamma_{1}:=\Gamma_{0}\cap {B_{\delta}(\bar{x})}\times {B_{D}(\bar{x})}$.
Since $\mu_{0}(B_{\delta}(\bar{x})) = \mu_{1}(B_{D}(\bar{x}))= 1$, it is easy to check that for any optimal plan $\bar{\pi}\in\Pi(\mu_{0},\mu_{1})$, it holds
\begin{align}\label{4.2}
\bar{\pi}(\Gamma_{1})=1.
\end{align}
Then we define
$$\Gamma:=\{(\gamma_{s},\gamma_{t})|\gamma\in \textmd{Geo}(X), 0\leq s \leq t\leq1,(\gamma_{0},\gamma_{1})\in \Gamma_{1}\}.$$
Define transport relation to be
$$R=\Gamma\cup\Gamma^{-1},$$
where $\Gamma^{-1}:=\{(x, y)\in X\times X\mid (y, x)\in \Gamma\}$.
Denote by $R(x)=\{y\mid (x,y)\in R\}$.
Define the associated transport set to be
$$\mathcal{T}_{e}:=P_{1}(R\setminus\{x = y\}),$$
and the set of branching points to be $A=A^{+}\cup A^{-}$, where
$$A_{+} := \{x\in \mathcal{T}_{e}\mid\exists z, w\in \mathcal{T}_{e}, (x, z), (x,w)\in \Gamma, (z,w)\notin R\},$$
$$A_{-} := \{x\in \mathcal{T}_{e}\mid\exists z, w\in \mathcal{T}_{e}, (x, z), (x,w)\in \Gamma^{-1}, (z,w)\notin R\},$$
and define the transport set without branching points to be
$$\mathcal{T}:= \mathcal{T}_{e}\setminus A.$$
One can check that the set $\mathcal{T}$ is Borel.
Making use of the essentially non-branching assumption, and the $\mathrm{MCP}(K,N)$ assumption, together with Theorem 1.1 in \cite{CaMo16}, we can follow the proof of Proposition 4.5 in \cite{Ca14-2} verbatim to obtain
\begin{align}\label{5.3}
\bar{m}(A) = 0.
\end{align}

In \cite{CaMo17}, it is proved that $R^{b} := R\cap(\mathcal{T}\times\mathcal{T})$ is an equivalence relation over $\mathcal{T}$, and for any $x\in\mathcal{T}$, the equivalence class $R^{b}(x)$ is isometric to an interval.
See Corollary 3.6 and Lemma 3.7 in \cite{CaMo17}.

There exists an $\mathcal{A}$-measurable map $\mathfrak{Q}:\mathcal{T}\rightarrow \mathcal{T}$ such that $(x, \mathfrak{Q}(x))\in R^{b}$ and $\mathfrak{Q}(x) = \mathfrak{Q}(y)$ whenever $(x, y)\in R^{b}$, and the quotient set $Q := \{x = \mathfrak{Q}(x)\}$ is $\mathcal{A}$-measurable.
See Lemma 3.8 in \cite{CaMo17}.
Here $\mathcal{A}$ denotes the $\sigma$-algebra generated by analytic sets.
Then we endow a Borel measure on $Q$ defined by
\begin{align}\label{4.30}
\mathfrak{q}=\mathfrak{Q}_{\#}(\bar{m}\mid _{\mathcal{T}}).
\end{align}

For $q\in Q$, we use $X_{q}$ to denote the equivalence class $R^{b}(q)$.
By construction, each $X_{q}$ is a geodesic, and it is part of a possibly longer geodesic whose two end points are contained in ${B_{\delta}(\bar{x})}$ and ${B_{D}(\bar{x})}$ respectively.
Hence by the triangle inequality, one can easily check that, for every $q\in Q$, the length of $X_{q}$ (denoted by $L_{q}$), is no larger than $D+\delta$,
and $X_{q}\subset {B_{D+2\delta}(\bar{x})}$.

\vspace*{10pt}

In the above we have introduced the notation when the localization method applied to $B_{D+2\delta}(\bar{x})$.
Now we give some important conclusions.

\begin{description}
  \item[(1)] $B_{D+2\delta}(\bar{x})$ can be written as the disjoint union of two sets $Z$ and $\mathcal{T}$ with $\mathcal{T}$ admitting a partition $\{X_{q}\}_{q\in Q}$; every $X_{q}$ is a geodesic in $(X, d)$ with $L_{q}\leq D+\delta$.
  \item[(2)] There exists a family of measures $\{\bar{m}_{q}\}_{q\in Q}\subset\mathcal{M}(X)$ such that, for  $\mathfrak{q}$-a.e. $q\in Q$, $\bar{m}_{q}$ is a probability measure and is concentrated on $X_{q}$; for every Borel set $C$, the map $q\mapsto \bar{m}_{q}(C)$ is $\mathfrak{q}$-measurable, and it holds
  \begin{align}\label{disint}
  \bar{m}\mid_{\mathcal{T}}(C)= \int_{Q}\bar{m}_{q}(C\cap\mathcal{T}) \mathfrak{q}(dq).
  \end{align}
  \item[(3)] For $\mathfrak{q}$-a.e. $q\in Q$, $\bar{m}_{q}=h_{q}\mathcal{H}^{1}\mid_{X_{q}}\ll \mathcal{H}^{1}\mid_{X_{q}}$, and $(X_{q}, d, \bar{m}_{q})$ is an $\mathrm{MCP}(K, N)$ space.
  \item[(4)] $f_{E}=0$ $\bar{m}$-a.e. in $Z$, where $f_{E}$ was defined in (\ref{5.1}).
  \item[(5)] For $\mathfrak{q}$-a.e. $q\in Q$, it holds
  \begin{align}\label{5.5}
  \int_{X_{q}}f_{E}\bar{m}_{q}=0.
  \end{align}
\end{description}

\vspace*{10pt}

Properties (1)-(5) are standard in the localization technique, as we briefly explain below.
In (1), the set $\mathcal{T}$, the map $\mathfrak{Q}:\mathcal{T}\rightarrow \mathcal{T}$, the section $Q$, the measure $\mathfrak{q}$ are obtained in the previous construction.
(2) is obtained by applying the disintegration theorem (see e.g. Section 452 in \cite{Frem-4} or Theorem A.7 in \cite{BC09}) to decompose $\bar{m}\mid_{\mathcal{T}}$ according to the quotient map $\mathcal{Q}$.
(3) can be obtained as in the proof of Theorem 9.5 in \cite{BC13}.
The proof of (4) and (5) can consult Step 2 and Step 3 in the proof of Theorem 5.1 in \cite{CaMo17I} respectively.
We remark that in these two part of proofs in \cite{CaMo17I}, the authors only use the fact (\ref{5.3}) and use some argument based on basic definitions in optimal transport.
Using (\ref{4.2}), we can slightly modify the proofs in \cite{CaMo17I} to obtain (4) and (5) in our setting.

\vspace*{10pt}

Besides (1)-(5), we need to supplement some properties which will be used in the proofs of Theorem \ref{thm-local-isoperi1} and Theorem \ref{thm-local-isoperi2}.

Following the proof of Theorem 7.10 in \cite{CaMi16}, we conclude that, for $\mathfrak{q}$-a.e. $q\in Q$, $\bar{X}_{q}$ coincide with $R(q)$, which is a geodesic whose two end points are contained in ${B_{\delta}(\bar{x})}$ and ${B_{D}(\bar{x})}$ respectively.
Thus we have
\begin{description}
  \item[(6)] For $\mathfrak{q}$-a.e. $q\in Q$, $X_{q}\cap B_{\delta}(\bar{x})\neq\emptyset$.
\end{description}

Combing (4) with the facts that $f_{E}(x)>0$ for $x\in E$ and $f_{E}(x)<0$ for $x\in B_{D}(\bar{x})\setminus E$, we have
$$\bar{m}(B_{D}(\bar{x})\cap Z) = 0.$$
Thus if $\delta$ is sufficiently small (depending on $N$, $K$ and $D$), then we have
\begin{align}\label{5.8}
&\bar{m}(\mathcal{T})\geq \bar{m}(B_{D}(\bar{x})) = \frac{m(B_{D}(\bar{x}))}{m(B_{D+2\delta}(\bar{x}))} \\
\geq & \frac{\mathrm{Vol}_{K,N}(D)}{\mathrm{Vol}_{K,N}(D+2\delta)}
=\frac{\int_{0}^{D}s_{K/(N-1)}(t)^{N-1}dt}{\int_{0}^{D+2\delta}s_{K/(N-1)}(t)^{N-1}dt} \nonumber\\
\geq & 1-C\delta,\nonumber
\end{align}
where $C$ is a positive constant depending on $N$, $K$ and $D$.
In conclusion, we have
\begin{description}
  \item[(7)] If $\delta$ is sufficiently small (depending on $N$, $K$ and $D$), then
  \begin{align}\label{5.9}
\mathfrak{q}(Q)=\bar{m}(\mathcal{T})\geq 1-C\delta.
\end{align}
\end{description}

\section{Proof of Theorem \ref{thm-local-isoperi1}}\label{sec4}

In this section, $C$ denotes some positive constant depending only on $N$, $K$, $D$, and it may vary in different lines.
Recall that $K=0$ or $-(N-1)$ and correspondingly, $\kappa=0$ or $-1$.
For every $L>0$, denote by
$$k_{L}=\int_{0}^{L}s_{\kappa}(t)dt.$$

By (\ref{5.5}) and (\ref{5.1}), we have
\begin{align}
0 =\int_{X_{q}}f_{E}\bar{m}_{q}=\bar{m}_{q}(E\cap X_{q}) -\frac{\bar{m}(E)}{\bar{m}(B_{D}(\bar{x}))}\bar{m}_{q}(B_{D}(\bar{x})), \qquad \text{for }\mathfrak{q}\text{-a.e. }q\in Q.
\end{align}

If $\rho$ is sufficiently small, then
$E\subset E^{\rho}\subset B_{D}(\bar{x})$.
Therefore, we have
\begin{align}\label{5.12}
\bar{m}^{+}(E)&=\liminf_{\rho\downarrow0}\frac{\bar{m}(E^{\rho})-\bar{m}(E)}{\rho}\\
&=\liminf_{\rho\downarrow0}\frac{\bar{m}(E^{\rho}\cap\mathcal{T})-\bar{m}(E\cap\mathcal{T})}{\rho}\nonumber\\
&=\liminf_{\rho\downarrow0}\int_{Q}\frac{\bar{m}_{q}(E^{\rho}\cap X_{q})-\bar{m}_{q}(E\cap X_{q})}{\rho}\mathfrak{q}(dq)\nonumber\\
&\geq\liminf_{\rho\downarrow0}\int_{Q}\frac{\bar{m}_{q}((E\cap X_{q})^{\rho})-\bar{m}_{q}(E\cap X_{q})}{\rho}\mathfrak{q}(dq)\nonumber\\
&\geq\int_{Q}\liminf_{\rho\downarrow0}\frac{\bar{m}_{q}((E\cap X_{q})^{\rho})-\bar{m}_{q}(E\cap X_{q})}{\rho}\mathfrak{q}(dq)\nonumber\\
&=\int_{Q}\bar{m}^{+}_{q}(E\cap X_{q})\mathfrak{q}(dq)\nonumber\\
&\geq \int_{Q}\tilde{\mathcal{I}}_{K,N,L_{q}}(\bar{m}_{q}(E\cap X_{q}))\mathfrak{q}(dq)\nonumber\\
&= \int_{Q}\tilde{\mathcal{I}}_{K,N,L_{q}}(\frac{\bar{m}_{q}(B_{D}(\bar{x}))}{\bar{m}(B_{D}(\bar{x}))}\bar{m}(E)) \mathfrak{q}(dq).\nonumber
\end{align}
In (\ref{5.12}), $(E\cap X_{q})^{\rho}:=\{x\in X_{q}\mid d(x, E\cap X_{q})<\rho\}$, $L_{q}$ denotes the length of $X_{q}$,
and we use Fatou's Lemma in the fifth line.

In the following, we assume $q\in Q$ satisfies all the properties in (2)-(6).

By (\ref{5.8}), $\bar{m}(B_{D}(\bar{x}))\geq 1-C\delta$, hence
\begin{align}\label{5.13}
\frac{\bar{m}_{q}(B_{D}(\bar{x}))}{\bar{m}(B_{D}(\bar{x}))}\leq\frac{1}{1-C\delta}\leq 1+C\delta.
\end{align}

If we view $X_{q}$ as a map of constant-speed parametrization $X_{q} : (0, L_{q}) \rightarrow X$ of the geodesic $X_{q}$, then, since $h_{q}$ is an $\mathrm{MCP}(K, N)$ density on $(0, L_{q})$ which integrates to $1$, by Lemma 2.4 in \cite{CaSa18}, it holds
\begin{align}\label{upperbound}
\sup_{x\in(0,L_{q})}h_{q}(x)\leq \frac{1}{L_{q}}
    \bigl(\int_{0}^{1}(\sigma_{K,N-1}^{(t)}(L_{q}))^{N-1}dt\bigr)^{-1}.
\end{align}

Suppose for some $q\in Q$ it holds $X_{q}^{-1}(B_{D+2\delta}(\bar{x})\setminus B_{D}(\bar{x}))\neq\emptyset$.
By property (6), $X_{q}$ intersects ${B_{\delta}(\bar{x})}$, hence $L_{q}\geq D-\delta> \frac{D}{2}$.
Also note that $L_{q}\leq D+\delta$ holds, hence by (\ref{upperbound}), for such $q$, we have
\begin{align}
\sup_{x\in(0,L_{q})}h_{q}(x)\leq C.
\end{align}

Thus
\begin{align}
\bar{m}_{q}(B_{D+2\delta}(\bar{x})\setminus B_{D}(\bar{x}))\leq\int_{D-\delta}^{L_{q}}h_{q}(x)dx\leq C\delta,
\end{align}
and hence
\begin{align}\label{5.16}
\bar{m}_{q}(B_{D}(\bar{x}))=\bar{m}_{q}(B_{D+2\delta}(\bar{x}))- \bar{m}_{q}(B_{D+2\delta}(\bar{x})\setminus B_{D}(\bar{x}))  \geq 1-C\delta.
\end{align}

If $X_{q}^{-1}(B_{D+2\delta}(\bar{x})\setminus B_{D}(\bar{x}))=\emptyset$, then (\ref{5.16}) still holds because in this case $\bar{m}_{q}(B_{D}(\bar{x}))=1$.

Thus we have
\begin{align}\label{5.17}
\frac{\bar{m}_{q}(B_{D}(\bar{x}))}{\bar{m}(B_{D}(\bar{x}))}\geq 1-C\delta.
\end{align}

\begin{rem}
If $K=N-1$, then by Lemma 2.4 in \cite{CaSa18}, we have
\begin{align}\label{rem5.1}
\sup_{x\in(0,L_{q})}h_{q}(x)\leq \frac{N}{L_{q}}.
\end{align}
Thus the above argument still holds and we still have (\ref{5.13}) and (\ref{5.17}).
\end{rem}

By (\ref{4.12}), we choose $\bar{v}$ depending on $N$, $K$ and $D$ such that
\begin{align}\label{5.10}
\tilde{\mathcal{I}}_{C,N,L}(v)= k_{L}^{-\frac{1}{N}}v^{\frac{N-1}{N}}(1+o(1))
\end{align}
holds for $v\in (0,2\bar{v})$ and $L\in[\frac{D}{2},2D]$.

Now we choose $\bar{\delta}$ sufficiently small (depending on $N$, $K$, $D$ and $\bar{f}$) so that,
\begin{align}
\sup_{r\in(0,\bar{\delta})\cap(0,r_{0})}\bar{f}(r)<\bar{v}
\end{align}
and for every $\delta\in(0,\bar{\delta})$,
\begin{align}\label{5.14}
\frac{k_{D+\delta}}{k_{D}}=\frac{\int_{0}^{D+\delta}s_{\kappa}(t)dt}{\int_{0}^{D}s_{\kappa}(t)dt}\leq 1+C\delta.
\end{align}

By (\ref{5.13}) and assumption (2),  we have:

\begin{align}\label{5.11}
\frac{\bar{m}_{q}(B_{D}(\bar{x}))}{\bar{m}(B_{D}(\bar{x}))}\bar{m}(E) \leq (1+C\delta)\bar{m}(E)<2\bar{v}.
\end{align}

By Lemma 3.9 in \cite{CaSa18}, the map $D\mapsto \tilde{\mathcal{I}}_{K,N,D}(v)$ is strictly decreasing because $K\leq0$.
Combining this fact with (\ref{5.17}), (\ref{5.10}) (\ref{5.14}) and (\ref{5.11}), we have

\begin{align}\label{5.18}
&\mathcal{I}_{K,N,L_{q}}(\frac{\bar{m}_{q}(B_{D}(\bar{x}))}{\bar{m}(B_{D}(\bar{x}))}\bar{m}(E))\geq \mathcal{I}_{K,N,D+\delta}(\frac{\bar{m}_{q}(B_{D}(\bar{x}))}{\bar{m}(B_{D}(\bar{x}))}\bar{m}(E)) \\
=& k_{D+\delta}^{-\frac{1}{N}}\biggl( \frac{\bar{m}_{q}(B_{D}(\bar{x}))}{\bar{m}(B_{D}(\bar{x}))}\bar{m}(E)\biggr)^{1-\frac{1}{N}} + o((\frac{\bar{m}_{q}(B_{D}(\bar{x}))}{\bar{m}(B_{D}(\bar{x}))}\bar{m}(E))^{1-\frac{1}{N}}) \nonumber\\
\geq &(1-C\delta)^{1-\frac{1}{N}}(1+C\delta)^{-\frac{1}{N}}k_{D}^{-\frac{1}{N}}\bar{m}(E)^{1-\frac{1}{N}} + o((\bar{m}(E))^{1-\frac{1}{N}})\nonumber \\
= &(1-\Psi(\delta))k_{D}^{-\frac{1}{N}}\bar{m}(E)^{1-\frac{1}{N}}  \nonumber
\end{align}

By (\ref{5.12}) (\ref{5.9}) and (\ref{5.18}), we have
\begin{align}
\bar{m}^{+}(E)\geq (1-C\delta)[(1-\Psi(\delta))k_{D}^{-\frac{1}{N}}\bar{m}(E)^{1-\frac{1}{N}}] =(1-\Psi(\delta))k_{D}^{-\frac{1}{N}}\bar{m}(E)^{1-\frac{1}{N}}.
\end{align}
Combined with assumption (1), we obtain
\begin{align}
m^{+}(E)&\geq m(B_{D+2\delta}(\bar{x}))^{\frac{1}{N}}(1-\Psi(\delta))k_{D}^{-\frac{1}{N}}m(E)^{1-\frac{1}{N}}\\
&\geq (N\omega_{N}k_{D})^{\frac{1}{N}}(1-\Psi(\delta))k_{D}^{-\frac{1}{N}}m(E)^{1-\frac{1}{N}}\nonumber\\
&\geq (1-\Psi(\delta))N^{\frac{1}{N}}\omega_{N}^{\frac{1}{N}}m(E)^{1-\frac{1}{N}}\nonumber.
\end{align}
This complete the proof.

\begin{rem}\label{rem5.2}
We consider a family of $1$-dimensional spaces $(X, d, m_{a}) =([0,D], |\cdot|, \tilde{h}_{a}\mathcal{L}^{1})$ (where $a>0$).
Here $\tilde{h}_{a}=\mathrm{Vol}_{K,N}(D)h_{a}$, with $h_{a}(x)=h^{a}_{K,N,D}(x)$ given in (\ref{4.3}), $K=-(N-1)$, $0$, or $N-1$.
Assume $a$ is sufficiently small (depending on $N$, $D$), then for $r\in[0, a]$, we have
\begin{align}
m_{a}([0,r])&=\mathrm{Vol}_{K,N}(D)\frac{f_{D}(a)}{s_{\kappa}(D-a)^{N-1}}\int_{0}^{r} s_{\kappa}(D-t)^{N-1}dt \\ &=N\omega_{N}k_{D} \frac{\frac{a^{N-1}}{k_{D}}(1+o(1))}{s_{\kappa}(D)^{N-1}+o(1)}\int_{0}^{r} (s_{\kappa}(D)^{N-1}+o(1))dt \nonumber\\
&=N\omega_{N}a^{N-1}r(1+o(1)) .\nonumber
\end{align}
Thus
$\limsup_{r\downarrow0}\frac{m_{a}([0,r])}{\omega_{N}r^{N}}=+\infty$ as $N>1$,
and
$\frac{1}{m_{a}([0,D])}m_{a}([0,r])\leq r$ provided $a$ is sufficiently small.
Choose $\bar{f}(r)=r$, $\bar{x}=0\in X$, then assumption (2) in Theorem \ref{thm-local-isoperi1} always holds for every $a>0$ sufficiently small.

For any $\epsilon>0$, by (\ref{4.12}),
\begin{align}\label{5.22}
\tilde{\mathcal{I}}_{K,N,D}(v) \geq (1-\epsilon)k_{D}^{-\frac{1}{N}}v^{\frac{N-1}{N}}
\end{align}
holds for every $v$ sufficiently small.
For every $v$ sufficiently small, we take $a=a_{D}(v)$ as in (\ref{4.11}), and then
take $(X, d, m_{a})$, $E=[0,a]$.
Note that $m_{a}(E)=\mathrm{Vol}_{K,N}v=N\omega_{N}k_{D}v$.
By Theorem \ref{4.1} and (\ref{5.22}), we have
\begin{align}
&m_{a}^{+}(E)=\mathrm{Vol}_{K,N}(D)h_{a}(a)=N\omega_{N}k_{D}\tilde{\mathcal{I}}_{K,N,D}(v)\\
\geq &(1-\epsilon)N\omega_{N} k_{D}^{\frac{N-1}{N}}v^{\frac{N-1}{N}} = (1-\epsilon)N^{\frac{1}{N}}\omega_{N}^{\frac{1}{N}} (m_{a}(E))^{\frac{N-1}{N}}.\nonumber
\end{align}
Thus the constant $N^{\frac{1}{N}}\omega_{N}^{\frac{1}{N}}$ in (\ref{1.3}) is sharp.

\end{rem}

\section{Proof of Theorem \ref{thm-local-isoperi2}}\label{sec5}

In this section, $K=N-1$ and $\kappa=1$.
For every $L\in(0,\pi)$, denote by
$$k_{L}=\int_{0}^{L}s_{\kappa}(t)dt.$$
$C$ denotes some positive constant depending only on $N$, $K$, $D$, and it may vary in different lines.

Let $Q_{1}:=\{q\in Q|L_{q}<\frac{D}{2}\}$, $Q_{2}:=Q\setminus Q_{1}$.
Denote by $A=\mathfrak{q}(Q_{1})$. 

Given $\delta>0$ sufficiently small, we choose $c=\delta^{-\frac{1}{2}}$.
Hence $(c+1)\delta<3\delta^{\frac{1}{2}}<\frac{D}{10}$.

By property (6), for $\mathfrak{q}$-a.e. $q\in\mathcal{Q}$, there exists $\check{q}\in X_{q}$ such that $\check{q}\in B_{\delta}(\bar{x})$.
Hence we have $B_{c\delta}(\check{q})\subset B_{(c+1)\delta}(\bar{x})$.
Thus by (\ref{disint}), we have
\begin{align}
&\bar{m}(B_{(c+1)\delta}(\bar{x}))\geq \int_{Q_{1}}\bar{m}_{q} (B_{c\delta}(\check{q}))\mathfrak{q}(dq) +\int_{Q_{2}}\bar{m}_{q} (B_{c\delta}(\check{q}))\mathfrak{q}(dq) \\
\geq &\int_{Q_{1}}\bar{m}_{q}(B_{\frac{D}{2}}(\check{q})) \frac{\mathrm{Vol}_{K,N}(c\delta)}{\mathrm{Vol}_{K,N}(\frac{D}{2})}\mathfrak{q}(dq) +\int_{Q_{2}}\bar{m}_{q} (B_{D+\delta}(\check{q})) \frac{\mathrm{Vol}_{K,N}(c\delta)}{\mathrm{Vol}_{K,N}(D+\delta)}\mathfrak{q}(dq) \nonumber\\
=&\int_{Q_{1}} \frac{\mathrm{Vol}_{K,N}(c\delta)}{\mathrm{Vol}_{K,N}(\frac{D}{2})}\mathfrak{q}(dq) +\int_{Q_{2}} \frac{\mathrm{Vol}_{K,N}(c\delta)}{\mathrm{Vol}_{K,N}(D+\delta)}\mathfrak{q}(dq), \nonumber
\end{align}
where in the second inequality, we apply Bishop-Gromov inequality to the $\mathrm{MCP}(K,N)$ spaces $(X_{q},|\cdot|,\bar{m}_{q})$, and use the fact that $B_{r}(\check{q})\cap X_{q}$ is identical to $\{y\in X_{q}\mid |y\check{q}|< r\}$ for $r>0$.

Then we have
\begin{align}\label{5.2}
&\frac{1}{m(B_{D+2\delta}(\bar{x}))} \frac{\mathrm{Vol}_{K,N}((c+1)\delta)}{\mathrm{Vol}_{K,N}(c\delta)} \frac{m(B_{(c+1)\delta}(\bar{x}))}{\mathrm{Vol}_{K,N}((c+1)\delta)} \\
\geq & \int_{Q_{1}} \frac{1}{\mathrm{Vol}_{K,N}(\frac{D}{2})} \mathfrak{q}(dq) +\int_{Q_{2}} \frac{1}{\mathrm{Vol}_{K,N}(D+\delta)} \mathfrak{q}(dq). \nonumber
\end{align}

If $\delta$ is sufficiently small, we have
\begin{align}\label{6.6}
\frac{\mathrm{Vol}_{K,N}((c+1)\delta)}{\mathrm{Vol}_{K,N}(c \delta)}\leq (\frac{c+2}{c})^{N} =1+2N\delta^{\frac{1}{2}}+o(\delta^{\frac{1}{2}})
\end{align}

By assumption (2) and Bishop-Gromov inequality, we have
\begin{align}
\frac{m(B_{(c+1)\delta}(\bar{x}))}{\mathrm{Vol}_{K,N}((c+1)\delta)}\leq 1+\eta.
\end{align}
Combining it with assumption (1) and (\ref{5.9}), (\ref{5.2}), (\ref{6.6}), we have

\begin{align}\label{6.8}
&\biggl(1+2N\delta^{\frac{1}{2}}+o(\delta^{\frac{1}{2}})\biggr) \frac{1+\eta}{(1-\eta)\mathrm{Vol}_{K,N}(D)} \\
\geq & \frac{A}{\mathrm{Vol}_{K,N}(\frac{D}{2})} +  \frac{1-C\delta-A}{\mathrm{Vol}_{K,N}(D+\delta)} \nonumber
\end{align}

Since $\delta$ and $\eta$ is sufficiently small, we have

\begin{align}\label{6.9}
\frac{\mathrm{Vol}_{K,N}(D)}{\mathrm{Vol}_{K,N}(D+\delta)}= \frac{k_{D}}{k_{D}+\int_{D}^{D+\delta}s_{\kappa}(t)dt} =1-\frac{1}{k_{D}}\delta s_{\kappa}(D)+o(\delta).
\end{align}

Denote by $h:=\frac{\mathrm{Vol}_{K,N}(D)}{\mathrm{Vol}_{K,N}(\frac{D}{2})}>1$.
By (\ref{6.8}) and (\ref{6.9}), we have

\begin{align}
&\biggl(1+2N\delta^{\frac{1}{2}}+o(\delta^{\frac{1}{2}})\biggr)\frac{1+\eta}{1-\eta} \\
\geq & Ah +  (1-C\delta-A)(1-\frac{1}{k_{D}}\delta s_{\kappa}(D)+o(\delta))  \nonumber\\
= & 1 -\frac{1}{k_{D}}\delta s_{\kappa}(D)+o(\delta)+ A(h-1+\frac{1}{k_{D}}\delta s_{\kappa}(D)+o(\delta)),\nonumber
\end{align}
and thus
\begin{align}
&A(h-1+\frac{1}{k_{D}}\delta s_{\kappa}(D)+o(\delta)) \\
\leq & \frac{1}{k_{D}}\delta s_{\kappa}(D)-1+o(\delta)+\biggl(1+2N\delta^{\frac{1}{2}}+o(\delta^{\frac{1}{2}})\biggr) \biggl(1+2\eta+o(\eta)\biggr) \nonumber \\
\leq & \frac{1}{k_{D}}\delta s_{\kappa}(D)+o(\delta)+2N\delta^{\frac{1}{2}}+o(\delta^{\frac{1}{2}})+2\eta \biggl(1+2N\delta^{\frac{1}{2}}+o(\delta^{\frac{1}{2}})\biggr) +o(\eta). \nonumber
\end{align}

If $\delta$ and $\eta$ are sufficiently small (depending on $N$, $D$), we have
\begin{align}
A(h-1)\leq 4N\delta^{\frac{1}{2}}+3\eta,
\end{align}
hence
\begin{align}\label{6.17}
A\leq \frac{4N\delta^{\frac{1}{2}}+3\eta}{h-1}.
\end{align}

By (\ref{5.12}), we have
\begin{align}\label{6.18}
\bar{m}^{+}(E)&\geq  \int_{Q}\tilde{\mathcal{I}}_{K,N,L_{q}}(\frac{\bar{m}_{q}(B_{D}(\bar{x}))}{\bar{m}(B_{D}(\bar{x}))}\bar{m}(E)) \mathfrak{q}(dq)\\
&\geq  \int_{Q_{2}}\tilde{\mathcal{I}}_{K,N,L_{q}}(\frac{\bar{m}_{q}(B_{D}(\bar{x}))}{\bar{m}(B_{D}(\bar{x}))}\bar{m}(E)) \mathfrak{q}(dq). \nonumber
\end{align}

By assumption (2) and Bishop-Gromov inequality, we have
\begin{align}
\bar{m}(E)\leq \bar{m}(B_{\delta}(\bar{x}))\leq \bigl(\limsup_{r\downarrow0}\frac{\bar{m}(B_{r}(\bar{x}))}{\mathrm{Vol}_{K,N}(r)} \bigr)\mathrm{Vol}_{K,N}(\delta)\leq (1+\eta)\mathrm{Vol}_{K,N}(\delta).
\end{align}
Recall that for $K=N-1$, (\ref{5.13}) and (\ref{5.17}) still hold for $\mathfrak{q}$-a.e. $q\in Q$ (see Remark \ref{rem5.1}).
Thus
\begin{align}
\frac{\bar{m}_{q}(B_{D}(\bar{x}))}{\bar{m}(B_{D}(\bar{x}))}\bar{m}(E) \leq (1+C\delta)(1+\eta)\mathrm{Vol}_{K,N}(\delta).
\end{align}

By (\ref{4.12}), we can choose $\bar{\delta}, \bar{\eta}>0$ sufficiently small (depending on $N$, $D$) such that, for any $v\in(0,\bar{v})$ with $\bar{v}=(1+C\bar{\delta})(1+\bar{\eta})\mathrm{Vol}_{K,N}(\bar{\delta})$,
\begin{align}
\tilde{\mathcal{I}}_{K,N,L}(v) =k_{L}^{-\frac{1}{N}}v^{\frac{N-1}{N}}(1+o(1))
\end{align}
holds for every $L\in[\frac{D}{2},D+\bar{\delta}]$.
(We assume $\bar{\delta}$ is sufficiently small so that $D+\bar{\delta}<\pi$.)

For $q\in Q_{2}$, we have $\frac{D}{2}\leq L_{q}\leq D+\delta$.
Thus for $\mathfrak{q}$-a.e. $q\in Q_{2}$ and for any $\delta\in(0,\bar{\delta})$, $\eta\in(0,\bar{\eta})$, we have

\begin{align}\label{6.20}
&\tilde{\mathcal{I}}_{K,N,L_{q}}(\frac{\bar{m}_{q}(B_{D}(\bar{x}))}{\bar{m}(B_{D}(\bar{x}))}\bar{m}(E))\\
= & k_{L_{q}}^{-\frac{1}{N}}\biggl( \frac{\bar{m}_{q}(B_{D}(\bar{x}))}{\bar{m}(B_{D}(\bar{x}))}\bar{m}(E)\biggr)^{1-\frac{1}{N}} + o((\frac{\bar{m}_{q}(B_{D}(\bar{x}))}{\bar{m}(B_{D}(\bar{x}))}\bar{m}(E))^{1-\frac{1}{N}}) \nonumber\\
\geq &(1-C\delta)^{1-\frac{1}{N}}k_{D+\delta}^{-\frac{1}{N}}\bar{m}(E)^{1-\frac{1}{N}} + o((\bar{m}(E))^{1-\frac{1}{N}})\nonumber \\
\geq &(1-C\delta)^{1-\frac{1}{N}}(1+C\delta)^{-\frac{1}{N}}k_{D}^{-\frac{1}{N}}\bar{m}(E)^{1-\frac{1}{N}} + o((\bar{m}(E))^{1-\frac{1}{N}}) \nonumber \\
\geq & (1- \Psi(\delta))k_{D}^{-\frac{1}{N}}\bar{m}(E)^{1-\frac{1}{N}}, \nonumber
\end{align}
where in the second inequality, we use the fact that $k_{L_{q}}\leq k_{D+\delta}$,
and in the third inequality, we use the fact that $k_{D+\delta}\leq (1+C\delta)k_{D}$ provided $\delta$ is sufficiently small.

Combining (\ref{5.9}), (\ref{6.17}), (\ref{6.18}) and (\ref{6.20}), we have
\begin{align}
\bar{m}^{+}(E)
\geq &(1-C\delta-\frac{4N\delta^{\frac{1}{2}}+3\eta}{h-1})(1- \Psi(\delta))k_{D}^{-\frac{1}{N}}\bar{m}(E)^{1-\frac{1}{N}}\\
\geq &(1- \Psi(\delta,\eta))k_{D}^{-\frac{1}{N}}\bar{m}(E)^{1-\frac{1}{N}}.\nonumber
\end{align}
Thus by assumption (1), we have
\begin{align}
m^{+}(E)&\geq m(B_{D+2\delta}(\bar{x}))^{\frac{1}{N}}(1-\Psi(\delta,\eta))k_{D}^{-\frac{1}{N}}m(E)^{1-\frac{1}{N}}\\
&\geq (N\omega_{N}k_{D})^{\frac{1}{N}} (1-\Psi(\delta,\eta))k_{D}^{-\frac{1}{N}}m(E)^{1-\frac{1}{N}}\nonumber\\
&\geq (1-\Psi(\delta,\eta))N^{\frac{1}{N}}\omega_{N}^{\frac{1}{N}}m(E)^{1-\frac{1}{N}}\nonumber.
\end{align}
The proof is completed.

\vspace*{10pt}



\begin{thebibliography}{99}


\bibitem{AG11} L. Ambrosio, N. Gigli, A user's guide to optimal transport. Modelling and Optimisation of Flows on Networks, Lecture Notes in Mathematics, 2062 (2011), Springer.




\bibitem{AGM16} L. Ambrosio, N. Gigli, S. Di Marino, Perimeter as relaxed Minkowski content in metric measure spaces. Nonlinear Analysis TMA, http://dx.doi.org/10.1016/j.na.2016.03.010.


\bibitem{BS10} K. Bacher, K.-T. Sturm, Localization and tensorization properties of the curvature-dimension condition for metric measure spaces. J. Funct. Anal., 259 (2010), 28-56.



\bibitem{BaRi18} D. Barilari, L. Rizzi, Sharp measure contraction property for generalized H-type Carnot groups. Commun. Contemp. Math., 20 (2018): 1750081, 24.





\bibitem{BaRi19} D. Barilari, L. Rizzi, Sub-Riemannian interpolation inequalities. Invent. math., 215 (2019), 977-1038.


\bibitem{BBG85} P. B\'{e}rard, G. Besson, S. Gallot. Sur une in\'{e}galit\'{e} isop\'{e}rim\'{e}trique qui g\'{e}n\'{e}ralise celle de Paul L\'{e}vy-Gromov. Invent. Math. 80 (1985), 295-308.


\bibitem{BC13} S. Bianchini, F. Cavalletti, The Monge problem for distance cost in geodesic spaces. Comm. Math. Phys., 318 (2013), 615-673.



\bibitem{BC09} S. Bianchini, L. Caravenna, On the extremality, uniqueness and optimality of transference plans. Bull. Inst. Math. Acad. Sin. (N.S.) 4 (2009), 353-454.


\bibitem{Ca18} F. Cavalletti, An overview of $L^{1}$-optimal transportation on metric measure spaces. in: N. Gigli (Ed.), Measure Theory in Non-Smooth Spaces. De Gruyter Open.






\bibitem{Ca14-2} F. Cavalletti, Monge problem in metric measure spaces with Riemannian curvature-dimension condition. Nonlinear Anal., 99 (2014), 136-151.


\bibitem{CaHu15} F. Cavalletti, M. Huesmann, Existence and uniqueness of optimal transport maps. Ann. Inst. H. Poincar\'{e} Anal. Non Lin\'{e}aire, 32 (2015), 1367-1377.



\bibitem{CaMi16} F. Cavalletti, E. Milman, The globalization theorem for the Curvature-Dimension condition, to appear in Invent. Math., arXiv.org/abs/1612.07623.


\bibitem{CaMo17I} F. Cavalletti, A. Mondino, Sharp and rigid isoperimetric inequalities in metric-measure spaces with lower Ricci curvature bounds, Invent. Math., 208 (2017), 803-849.






\bibitem{CaMo16} F. Cavalletti, A. Mondino, Optimal maps in essentially non-branching spaces, Commun. Contemp. Math., 19 (2017), pp. 1750007, 27.

\bibitem{CaMo18} F. Cavalletti, A. Mondino, New formulas for the Laplacian of distance functions and applications, to appear in Anal. PDE, arXiv.org/abs/1803.09687.



\bibitem{CaMo18I} F. Cavalletti, A. Mondino, Isoperimetric inequalities for finite perimeter sets under lower Ricci curvature bounds, Rend. Lincei Mat. Appl., 29 (2018), 413-430.




\bibitem{CaMo17} F. Cavalletti, A. Mondino, Almost Euclidean Isoperimetric Inequalities in spaces satisfying local Ricci curvature lower bounds, Int. Math. Res. Not. IMRN, 2020 (2020), Issue 5, 1481-1510.



\bibitem{CaSa18} F. Cavalletti, F. Santarcangelo, Isoperimetric inequality under measure contraction property, J. Funct. Anal., 277 (2019), 2893-2917.



\bibitem{Cr82} C. B. Croke, An eigenvalue pinching theorem. Invent. Math., 68 (1982), 253-256.


\bibitem{EM13} M. Eichmair, J. Metzger, Unique isoperimetric foliations of asymptotically flat manifolds in all dimensions , Invent. Math., 194 (2013), 591-630.

\bibitem{Frem-4} D. H. Fremlin, Measure Theory, volume 4, Torres Fremlin (2002).





\bibitem{Gro07} M. Gromov, Metric structures for Riemannian and non Riemannian spaces. Modern Birkh\"{a}user Classics, (2007).

\bibitem{GrMi87} M. Gromov, V. D. Milman, Generalization of the spherical isoperimetric inequality to uniformly convex Banach spaces. Compos. Math., 62 (1987), 263-282.



\bibitem{H19} B.-X. Han, Sharp $p$-Poincar\'{e} inequality under Measure Contraction Property. Manuscripta Math., 162 (2020), 457-471.


\bibitem{HaMi18} B.-X. Han, E. Milman, Sharp Poincar\'{e} inequalities under Measure Contraction Property. To appear in Ann. Sc. Norm. Super. Pisa CI. Sci.



\bibitem{Jui09} N. Juillet, Geometric inequalities and generalized Ricci bounds in the Heisenberg group. Int. Math. Res. Not. IMRN, 13 (2009), 2347-2373.





\bibitem{KLS87} R. Kannan, L. Lov\'{a}sz, M. Simonovits, Isoperimetric problems for convex bodies and a localization lemma. Discrete Comput. Geom., 13 (1995), 541-559.



\bibitem{Kla17} B. Klartag, Needle decomposition in Riemannian geometry. Mem. Amer. Math. Soc., 249 (2017).


\bibitem{LV09} J. Lott, C. Villani, Ricci curvature for metric-measure spaces via optimal transport. Ann. of Math. (2), 169 (2009), 903-991.





\bibitem{LLZ16} P. W. Y. Lee, C. Li, I. Zelenko, Ricci curvature type lower bounds for sub-Riemannian structures on Sasakian manifolds. Discrete Contin. Dyn. Syst., 36 (2016), 303-321.


\bibitem{Mil15} E. Milman, Sharp isoperimetric inequalities and model spaces for curvature-dimension diameter condition. J. Eur. Math. Soc., 17 (2015), 1041-1078.



\bibitem{Oh07} S. Ohta, On the measure contraction property of metric measure spaces. Comment. Math. Helv., 82 (2007), 805-828.





\bibitem{PaWe60} L. E. Payne, H. F. Weinberger, An optimal Poincar\'{e} inequality for convex domains. Arch. Ration. Mech. Anal, 5 (1960), 286-292.



\bibitem{Per02} G. Perelman, The entropy formula for the Ricci flow and its geometric applications. arXiv:math/0211159v1, (2002).




\bibitem{St06I} K. T. Sturm, On the geometry of metric measure spaces I, Acta Math. 196 (2006), 65-131.

\bibitem{St06II} K. T. Sturm, On the geometry of metric measure spaces II, Acta Math. 196 (2006), 133-177.

\bibitem{Vi09} C. Villani, Optimal transport, Old and new. Grundlehren der Mathematischen Wissenschaften, Springer-Verlag, Berlin, vol. 338 (2009).








\end{thebibliography}
\end{document}